# Optimal Control of Chromate Removal via Enhanced Modeling using the Method of Moments

Fred Ghanem and Kirti M. Yenkie*

Department of Chemical Engineering, Henry M. Rowan College of Engineering,
Rowan University, 600 North Campus Drive, Glassboro, NJ08028, USA

(Tel: 856-256-5375; e-mail: yenkie@rowan.edu*; ghanem@students.rowan.edu)

**Abstract:**

Single-use anion-exchange resins can reduce hazardous chromates to safe levels in drinking water. However, since most process control strategies monitor effluent concentrations, detection of any chromate leakage leads to premature resin replacement. Furthermore, variations in the inlet chromate concentration and other process conditions make process control a challenging step. In this work, we capture the uncertainty of the process conditions by applying the Ito process of Brownian motion with a drift into a stochastic optimal control strategy. The ion exchange process is modeled using the method of moments which helps capture the process dynamics, later formulated into mathematical objectives representing desired chromate removal. We then solved our developed models as an optimal control problem via Pontryagin's maximum principle. The objectives enabled a successful control via flow rate adjustments leading to higher chromate extraction. Such an approach maximizes the capacity of the resin and column efficiency to remove toxic compounds from water while capturing deviations in the process conditions.

**Keywords:** Ion exchange modeling; Method of moments; Optimal control; Pontryagin's maximum principle; maximized resin usage; Stochastic Ito process

# 1. Introduction

Optimal control of single-use ion exchange resins is necessary to avoid leakage of hazardous chemicals such as chromates that contaminate our drinking water supply. While we can improve the model accuracy to capture the process characteristics, running the system at a constant flow rate leads to inefficient use of the ion exchange resin, therefore increasing the cost and disposal of solid resin waste. In this body of work, optimality is achieved by deviating the flow rate with time, a unique approach to maximizing the removal of chromate.

Chromium-based compounds, such as chromates and dichromates, are widely used in the plating industry to protect metals against corrosion and as nanocoating in aeronautical applications (Ferreira, Zheludkevich, and Tedim 2011). However, chromium is a highly oxidizing metal that can be carcinogenic to human health at even lower concentrations of parts per billion (ppb) (Costa and Klein 2006). One of the most common types of resins used to remove chromates is the strong base anion exchanger (Brito F et al. 1997).

The ion exchange process is the equivalent exchange of ions of a given charge (either positive cations or negative anions) in a solution with ions of the same charge released from solid material or the ion exchanger (de Dardel and Arden 2008). Standard ion exchange resins are polymeric materials used in a wide range of applications, such as hardness removal from drinking water (de Dardel and Arden 2008) and the extraction of toxic endotoxins from plasma in medical dialysis devices (Nagaki et al. 1991). Ion exchange resins are also used as the last step in producing semiconductor-grade water with less than 1 part per trillion (ppt) of ions in the effluent, achieving the low conductivity needed for rinsing microchips during manufacturing (Hutcheson 2006). As a result, ion-exchange resins are important to achieve high-purity effluent and the removal of hazardous materials from drinking water such as lead (Biswas and Mishra 2015), arsenic (Karakurt, Pehlivan, and Karakurt 2019), nitrates (Nur et al. 2015), and chromates (Lin and Kiang 2003) down to parts per billion level (Kabir 2008). Strong base anion exchangers have gained tremendous applications, especially in the tannery wastewater treatment (Kabir 2008) and waste brine cleanup (Li et al. 2016) of chromates. In this work, we evaluate the performance of a standard gel-type strong anion exchanger, Purolite A600E (Purolite Website 2021), for removing chromates from water.

Since any leakage of chromates is to be avoided, feedforward control is essential to minimize such occurrences (Corder and Lee 1986). But the use of such a control strategy requires good system models to predict the output and therefore modify the process before any chromate leakage. However, inaccurate models lead to erroneous predictions that could lead to early breakthroughs of toxic substances. Moreover, most process models have empirically calculated parameters that deviate from process changes, such as flow rates (Recepoğlu et al. 2018) and inlet concentrations (Charola et al. 2018). Hence, improved predictive models are needed to capture such sensitive behavior to avoid early system failure. Modeling challenging processes via used of minimization of the error difference between the measured and simulated data are necessary in many applications (Han et al. 2023). Previous work was published by the authors to increase the accuracy of such models in chromate extraction (Ghanem, Jerpoth, and Yenkie 2022). Since all purification systems depend on successful process control, predicting the dynamics and subsequently implementing a better control strategy will help to design efficient and smaller-scale ion exchange systems.

However, the numerical methods used for obtaining such prediction are based on complex partial differential equations which turn out to be difficult to implement as a model (LeVan, Carta, and Yon 1999). Therefore, many simpler kinetic models are currently being used for process prediction, such as Wolborska (Hamdaoui 2009), Yoon-Nelson (Millar et al. 2015), Clark (Hamdaoui 2009), and Thomas models (Kalaruban et al. 2016), most being sigmoidal function representations.

These models also depend on empirical parameters that are not fixed when system characteristics, such as flow rates and inlet concentrations, fluctuate (Recepoğlu et al. 2018). In this work, we use the Thomas model to represent the chromate removal process due to its accuracy when compared to other models, as described in detail in prior work (Ghanem, Jerpoth, and Yenkie 2022).

With the adoption of such predictive models, the method of moments, initially set up by Goltz and Roberts in 1987 for chromatographic separations (Goltz and Roberts 1987; Staby et al. 2007), is used to interpret the distribution curves of the output concentration. In addition, the method of moments has been implemented in other fields to represent particulate formations in crystallization (Kirti M. Yenkie and Diwekar 2012). Typically, the sigmoidal functions, presented as cumulative curves, need to be transformed into distribution curves before applying the method of moments. Specific temporal moments, such as the normalized first moment (residence time), the central second moment (variance), and the central third moment (skewness), can then be calculated when dealing with step input processes (Yu, Warrick, and Conklin 1999) and as described in the method section of this work.

To ensure the safety of water for consumption, dynamic optimal control strategies (Boscain and Piccoli 2005) gained interest when certain objectives such as maximum removal of hazardous chemicals or longer process run-time are desired. With the moments' equations utilized to describe the ion exchange system, Pontryagin's maximum principle (Artstein 2011) can be applied via the Hamiltonian system (Hamill 2018) which is a canonical transformation of a time-dependent system, as described by Diwekar (Diwekar 2008). Pontryagin's maximum principle (Harmand et al. 2019), generally used in optimal control strategy to regulate a dynamic system, is defined by the following three characteristics:

1) The dynamic system is the center of the control strategy, represented by the ion exchange resin via the method of moments, which are used as coordinates in the Hamiltonian system.
2) The control law is a time-dependent variable that can be adjusted to create changes to the system. In this report, the flow rate is an example of such a controlling parameter. This control variable can also be constrained by the system pump limitations, such as maximum or minimum flow rate.
3) The optimality criterion is the objective function needed to maximize the performance of the ion exchange system, such as the highest chromate removal capacity or the longest production time. The objective function uses the Hamiltonian function by minimizing its deviations.

Since the ion exchange resin, used to remove the hazardous chromates, is rarely regenerated due to the high toxicity of such chemicals, it is economically and environmentally important to minimize the volume of the resin needed for extraction. Currently, most processes built today are based on the worst-case scenarios of inlet chromate levels leading to system overdesign resulting in inefficient processes and more resin disposal. However, with the successful application of optimal control via the method of moments and Pontryagin's maximum principle, smaller systems can be used efficiently, leading to less solid resin waste. While the use of simulated moving bed systems (Sharma 2022) can lead to the same advantages, the investment in such equipment is higher than applying flow rate control onto a standard packed bed system, as we demonstrate in this body of work.

**2. Methods, Theory, and Calculations:**

The method and theory of this work are presented in eight subsequent sections. Section 2.1 introduces the predictive model sigmoidal function, section 2.2 sets up the method of moments, section 2.3 formulates the Hamiltonian function, section 2.4 expands into the deterministic optimal control using Pontryagin's maximum principle, section 2.5 proposes possible objective functions, section 2.6 captures uncertainties via Ito process, section 2.7 expands the use of the Ito process into a stochastic control strategy, and section 2.8 summarizes the entire approach by running both deterministic and stochastic optimization programs.

**2.1. Setting up the empirical model**: The experiments consisted of periodic measurements of the effluent chromate concentration when a constant inlet chromate feed is introduced. The experimental data were collected from the work of Xue Li and coworkers (Li 2016). As mentioned in the introduction, the Thomas model (Mustafa and Ebrahim 2010), a simple representation of the purification system, was used due to its higher accuracy and simple representation of the system characteristics (seen in equation (1)).

$$\Psi = \frac{C}{C_0} = \frac{1}{1+\exp\left(\frac{K_T q_m V}{Q} - K_T C_0 t\right)} \tag{1}$$

where $\psi$ – the normalized outlet chromate concentration to the inlet chromate concentration = $C/C_0$, $C_0$ - initial concentration of the chromate in the feed to the anion exchanger, C - outlet concentration predicted by the model, V - the volume of resin in the column, Q - flow rate of the feed, t – time elapsed since feed introduction, $K_T$ and $q_m$ – Thomas model parameters calculated by minimizing the sum of square error between the empirical and experimental results.

With changes in the experimental chromate concentrations as well as the flow rates to the ion exchange resin column, the empirical model was optimized in previous work and summarized below (Ghanem, Jerpoth, and Yenkie 2022). In that work, it was demonstrated that fixed parameters cannot capture accurately fluctuations in inlet chromate concentration and flow rates.

Therefore, the $q_m$ parameter was fixed while allowing the $K_T$ parameter to change linearly with chromate concentrations and contact times. The values for such parameters were calculated by regression analysis to be $q_m$ = fixed value = 0.254 kg Cr/m³ and $K_T$ having a relationship with contact time and inlet chromate concentration, $C_0$, as listed below:

$$K_T = \alpha\, CT + \beta\, C_0 + \gamma = -264\, CT + 10.45\, C_0 + 1247 \qquad (2)$$

with $K_T$ having units of l/gram Cr-hr, CT – contact time ($CT = V/Q$), having units of min, and $C_0$ having units of ppb ($\alpha$, $\beta$, and $\gamma$ being constants).

From equation (1), the following items are examples of variations that can influence the effluent concentration:

**i- variations in the inlet concentrations:** the concentration could range from 1% to 10% since chromate concentrations are in the parts of billion (Balan, Volf, and Bilba 2013). Therefore a 1-ppb variation, while small, is considered a dramatic variation that can lead to a premature ion exchange failure. Therefore, the inlet chromate concentration is the most important system characteristic that can affect the performance of the system (Ghanem, Jerpoth, and Yenkie 2022).

**ii- variations with the predictive model parameters, $K_T$ and $q_m$:** In previous work (Ghanem, Jerpoth, and Yenkie 2022), the accuracy of the parameters to reflect the system conditions were seen to fall within a 30% window, which needs to be captured in our stochastic process representation.

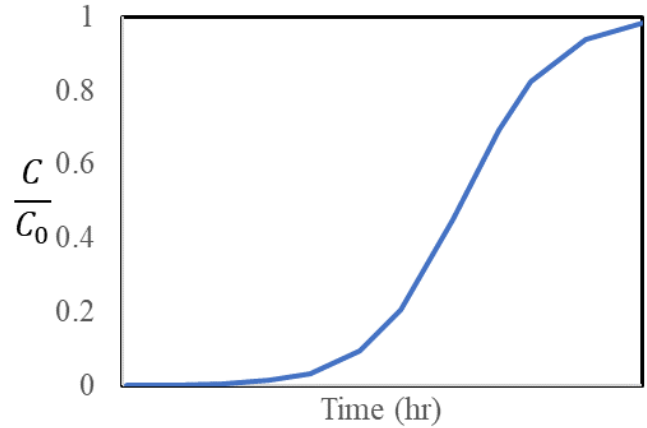

Figure 1. The sigmoidal curve of the concentration ratio, $C/C_0$, is based on a constant input concentration of $C_0$ with time.

Since the sigmoidal Thomas model is represented by a cumulative curve, as seen in Figure 1, it is necessary to transform such a model into a distribution curve, typically associated with chromatographic applications (Staby et al. 2007). That is accomplished by taking the derivative of the Thomas model from equation (1).

$$\frac{d\psi}{dt} = \frac{K_T C_0 \exp(K_T q_m CT - K_T C_0 t)}{(1+\exp(K_T q_m CT - K_T C_0 t))^2} = K_T C_0 (\psi - \psi^2) \qquad (3)$$

where $d\psi/dt$ becomes the new predictive model representation applied with the method of moments (Yu, Warrick, and Conklin 1999). The transformation of a sigmoidal or cumulative curve to a distribution curve is shown in Figure 2. As seen in the figure, the inflection point of the sigmoidal Thomas model representation becomes the peak point for the distribution curve while the sharpness of the increase around the inflection point becomes the width of the peak. Therefore, the derivative of the step input can be treated as a chromatographic separation where the peak time is related to the residence time and the peak width is related to the standard deviation.

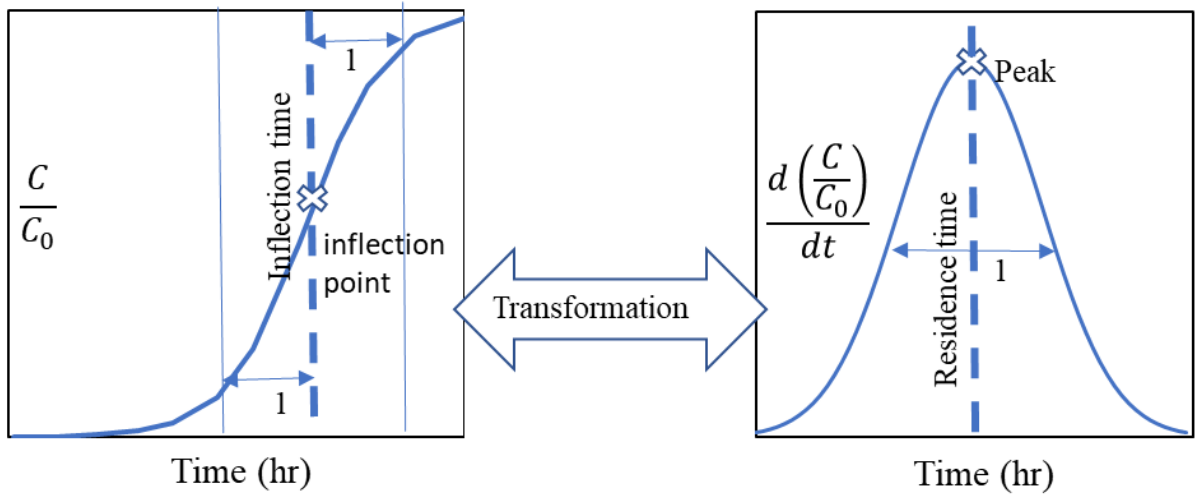

Figure 2. Temporal moments relationship to modeling the effluent concentration for a chromatographic pulse injection (fig, a) and constant input injection (fig. b) with "1" showing the variance related to the second central moment and peak/inflection point showing the first normalized moment (fig. a and fig. b).

**2.2. Using the method of moments**: As the ion exchange application is a fixed bed process, the concentration will change as it travels the length of the column with time. Therefore, the concentration is a function of bed depth and time. Partial differential equations (PDEs) are typically used to display the purification system as seen below.

$$\frac{\partial C}{\partial t} = D\frac{\partial^2 C}{\partial x^2} - v\frac{\partial C}{\partial x} - \left(\frac{1-\varepsilon}{\varepsilon}\right)\frac{\partial q}{\partial t} \quad \text{and} \quad \frac{\partial q}{\partial t} = \alpha_q(C - q) \qquad (4)$$

where C – Concentration in the effluent, D – dispersion coefficient, $v$ – intrinsic velocity in the column, $\varepsilon$ – void volume in a packed column, t – time, x – resin bed depth, $\alpha_q$ - first-order mass transfer rate constant, q – concentration accumulating on the resin. Rather than using such PDEs to represent the system, ordinary differential equations (ODEs), in the form of moments, are utilized as univariate time-dependent equations in this body of work.

The generalized method of moments is defined by the following equation for a chromatographic separation (Goltz and Roberts 1987):

$$m_i(\lambda) = \int_0^\infty t^i C(\lambda, t)\, dt \qquad (5)$$

with $m_i(t)$ - ith dimension temporal moment that is a function of $\lambda$, $\lambda$ - normalized column depth = x/L where $\psi$ value is predicted, $C(\lambda, t)$ –outlet concentration as a function of $\lambda$ and time, x – column depth where $\psi$ value is predicted, L – total resin depth, and t - time. Since we are predicting the concentration at the column outlet where $\lambda = 1$, the value of the moments in equation (4) will be constants (integrating between 0 to $\infty$). Therefore, the truncated temporal moments are used instead (Luo, Cirpka, and Kitanidis 2006) which measures the distribution between time = 0 and time = t rather than between 0 and $\infty$:

$$m_i(t) = \int_0^t t^i C(t) \, dt \qquad (6)$$

As previously mentioned, the moments, summarized by equation (6) represent a chromatographic process where the output is a distribution curve. Since the sigmoidal Thomas model, which is a step input process model imitating a population curve rather than a distribution curve as shown in Figure 2, can be rewritten using its derivative to display a chromatographic separation process as shown in equation (3):

$$m_i(t) = \int_0^t t^i \frac{d\psi}{dt} \, dt \qquad (7)$$

In chromatographic separation with a distribution curve output, four moments can be used to represent such a purification step: the zeroth moment, the first moment (reflecting the mean residence time of the distribution curve peak), the central second moment (the variance of the peak), and the central third moment (skewness of the peak) (Fogler 2016). Equation (7) is then used to calculate the truncated moment expressions for the four moments in the following equations (Jawitz 2004):

$$\text{Zeroth moment} = m_0(t) = \mu_0(t) = \int_0^t \frac{d\psi}{dt} \, dt = \psi(t) = \text{Thomas model} \qquad (8)$$

$$\text{First normalized moment} = \mu_1(t) = \frac{\int_0^t t \frac{d\psi}{dt} \, dt}{\int_0^t \frac{d\psi}{dt} \, dt} = \frac{\int_0^t t \frac{d\psi}{dt} \, dt}{\mu_0(t)} = t_m(t) \qquad (9)$$

$$\text{Second central normalized moment} = \mu_{2c}(t) = \frac{\int_0^t (t-\mu_1)^2 \frac{d\psi}{dt} dt}{\mu_0(t)} = \sigma^2(t) \qquad (10)$$

$$\text{Third central normalized moment} = \mu_{3c}(t) = \frac{1}{\mu_{2c}^{1.5}} \frac{\int_0^t (t-\mu_1)^3 \frac{d\psi}{dt} dt}{\mu_0(t)} = s^3(t) \qquad (11)$$

where $t_m$ – residence time for the maximum peak of $d\psi/dt$, $\sigma^2$ – variance or square of the standard deviation for $d\psi/dt$, and $s^3$ – skewness cubed for $d\psi/dt$. These four moments, typically measured via the output of the system, capture the system characteristics of the purification process. But as they are four separate equations, combining them into a single equation, such as a Lagrange expansion (Hamill 2018) with dynamic coordinates rather than cartesian ones, helps to modify the process into a simple representation.

**2.3. The Hamiltonian – a Lagrangian expansion:** The Hamiltonian system is a mathematical formulation that describes the dynamics of a physical system requiring some control strategies using the system coordinates rather than the standard cartesian representation (Hamill 2018). Hence, the Hamiltonian can be used to transform the temporal moments from four independent equations into coordinates of the new purification system. The Hamiltonian system, now defining the ion exchange process, can be represented by the following generalized equation (Diwekar 2008):

$$H = \sum_{i=1}^{n} z_i \frac{dy_i}{dt} = \sum_{i=1}^{n} z_i f(y_i, t, \chi) \qquad (12)$$

where H – Hamiltonian, $y_i$ – system coordinate representing the state variable or moment defined as $dy_i/dt = dH/dz_i$ (Chiang 1992), $z_i$ – system coordinate representing the adjoint variable defined as $dz_i/dt = -dH/dy_i$ (Chiang 1992), $\chi$ - variable to control in the process (such as the flow rate), and n – number of state equations being used in the optimization process. Therefore, the Hamiltonian represents the pathway the ion resin process will take based on its coordinates, $y_i$, which is defined by the moments, and $z_i$, which restricts the pathway of the system.

As mentioned, the state variables are the n equations that describe the ion exchange process. In this work, four state variables are defined by equation (8) for $y_1$, (9) for $y_2$, (10) for $y_3$, and (11) for $y_4$. As seen in the Hamiltonian in equation (12), the derivatives of the deterministic state variables are needed since the objective is to minimize system output deviations via pathway restriction.

$$\frac{dy_1}{dt} = \frac{d\psi}{dt} = K_T C_0 (y_1 - y_1^2) \tag{13}$$

$$\frac{dy_2}{dt} = \frac{(t - y_2)}{y_1} \frac{dy_1}{dt} \tag{14}$$

$$\frac{dy_3}{dt} = \frac{[(t - y_2)^2 - y_3]}{y_1} \frac{dy_1}{dt} \tag{15}$$

$$\frac{dy_4}{dt} = \frac{(t - y_1)^3}{y_1 y_3^{1.5}} \frac{dy_1}{dt} - y_4 \left( \frac{1}{y_1} \frac{dy_1}{dt} + \frac{1.5}{y_3} \frac{dy_3}{dt} \right) \tag{16}$$

As seen in equations (13) to (16), the derivatives of the state variables are written as functions of the state variables, $y_i$, enabling easier calculations of the adjoint variables in the next steps. The adjoint variables are used to counter the variation of the state variables, restricting the path of the Hamiltonian function, and are defined by the following generalized equation (Kirti M. Yenkie and Diwekar 2012):

$$\frac{dz_i}{dt} = -\sum_{k=1}^{n} z_k \frac{df(y_k, t, \chi)}{dy_i} \tag{17}$$

Equation (17) shows that the adjoint variables are also introduced through their derivatives, an important characteristic since derivatives describe a motion or change to achieve a desired result. Therefore equation (17) is used in conjunction with equations (13) to (16) to derive the 4 adjoint variables:

$$\frac{dz_1}{dt} = -z_1 \frac{d\left(\frac{dy_1}{dt}\right)}{dy_1} - z_2 \frac{d\left(\frac{dy_2}{dt}\right)}{dy_1} - z_3 \frac{d\left(\frac{dy_3}{dt}\right)}{dy_1} - z_4 \frac{d\left(\frac{dy_4}{dt}\right)}{dy_1} \tag{18}$$

$$\frac{dz_2}{dt} = -z_2 \frac{d\left(\frac{dy_2}{dt}\right)}{dy_2} - z_3 \frac{d\left(\frac{dy_3}{dt}\right)}{dy_2} - z_4 \frac{d\left(\frac{dy_4}{dt}\right)}{dy_2} \tag{19}$$

$$\frac{dz_3}{dt} = -z_3 \frac{d\left(\frac{dy_3}{dt}\right)}{dy_3} - z_4 \frac{d\left(\frac{dy_4}{dt}\right)}{dy_3} \tag{20}$$

$$\frac{dz_4}{dt} = -z_4 \frac{d\left(\frac{dy_4}{dt}\right)}{dy_4} \tag{21}$$

Equations (18) to (21) show that changes in the adjoint variables with time are running opposite to the shifts of the state variables, essential for a controlling parameter to minimize Hamiltonian system changes.

**2.4. Objective function:** Since the main objective of the ion exchange process is to maximize the removal of chromates, Figure 3 shows that an increase in the normalized first moment, $t_m$, results in more chromate removal, while a decrease in the second central moment, $\sigma^2$, results with lower leakage. Therefore, the objective for the best use of the ion exchange process would be to maximize $t_m$ and minimize the variance $\sigma^2$ simultaneously. If the weight is given equally to both moments, one possible objective function to maximize is the residence time with minimum variance shown below:

$$\text{Objective function} = t_m - \sigma \tag{22}$$

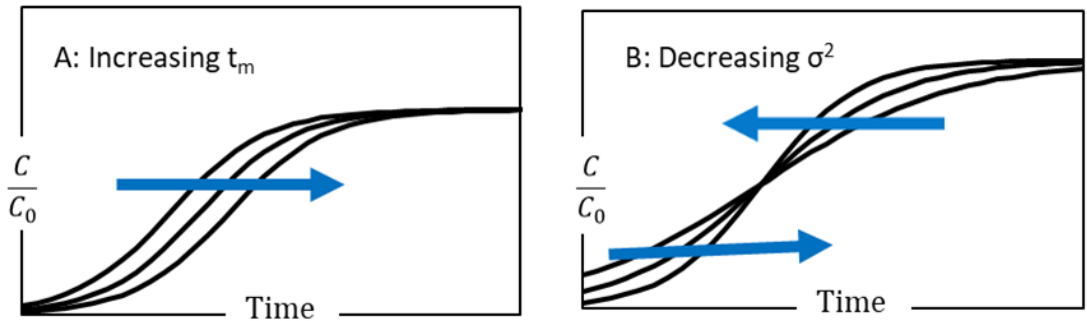

Figure 3. Effect of increasing the first normalized moment or residence time (A) compared to decreasing the second central moment or variance (B).

In chromatography, better purification is achieved when either the theoretical plate number (TPN) is maximized or the height equivalent to a theoretical plate (HETP) is minimized. HETP and TPN are related to the first and second moments via the following equations (Carta and Jungbauer 2010):

$$HETP = \left(\frac{\sigma^2}{t_m^2}\right) L = \frac{L}{TPN} \qquad (23)$$

where L = total height of the ion exchange resin column, and $TPN = t_m^2/\sigma^2$.

Equation (23) also demonstrates that TPN is maximized when the residence time is maximized at the same time as the variance is minimized, a similar objective demonstrated in the previous equation (22). But the TPN is typically maximized based on running the process to its completion when $C = C_0$. Since we cannot allow any large breakthrough of chromates, we aim to maximize the extraction until a certain breakthrough is achieved or the time has passed.

Table 1. The calculated dynamic binding capacity of the anion exchange resin from experimental data at different objective functions.

| Calculated from n experimental data points | Contact Time = 0.5 min | Contact Time = 1.5 min |
|---|---|---|
| $t_m$ = Residence time = $\sum_1^n t \frac{d\psi}{dt} dt$ | 108 hrs | 318 hrs |
| $\sigma^2$ = Variance = $\sum_1^n (t - t_m)^2 \frac{d\psi}{dt} dt$ | 3386 hrs² | 7150 hrs² |
| Objective Time 1 = $t_m - \sigma$ | 50 hrs | 233 hrs |
| C/C₀ at Objective Time 1 | 0.18 | 0.14 |
| Objective Time 2 = $t_m - \sigma/2$ | 78 hrs | 274 hrs |
| C/C₀ at Objective Time 2 | 0.32 | 0.28 |
| Objective Time 3 = $t_m + \sigma/2$ | 137 hrs | 359 hrs |
| C/C₀ at Objective Time 3 | 0.69 | 0.70 |

In this work, we chose an objective function that combined three moments: maximizing chromate removal by the time objective presented by equation (22), reflecting a breakthrough capacity of 14 to 18% as shown in Table 1. Table 1 shows the calculated concentration ratio at 3 objective times for 2 different contact times (short and long contact times = boundaries of our current chromate system). We picked the objective time of $t_m - \sigma$ to control the concentration ratio around 0.14 to 0.18 keeping the chromate levels under control.

The new objective noted as $J_{max}$ is seen below:

$$J_{max}(Q) = (C_0 - C)\, Q\, (t_m - \sigma) = C_0\bigl(1 - y_1(Q)\bigr)\, Q\, \bigl(y_2(Q) - \sqrt{y_3(Q)}\bigr) \quad (24)$$

where Q – flow rate, C – effluent concentration, and $C_0$ – inlet concentration. As mentioned in Table 1, equation (24) captures the amount of chromate extracted at $t_m - \sigma$.

**2.5. Deterministic optimal control and boundary conditions**: Figure 4 summarizes our approach for applying control strategies to the ion exchange resin system. It reflects our use of Pontryagin's maximum principle, described in the introduction, to provide optimal control for the ion exchange system via the temporal moments rather than the partial differential equations describing the system. This approach reflects the work done by Diwekar-Yenkie in crystallization (Kirti M. Yenkie and Diwekar 2012) and Diwekar-Benavides in biodiesel production (Benavides and Diwekar 2013).

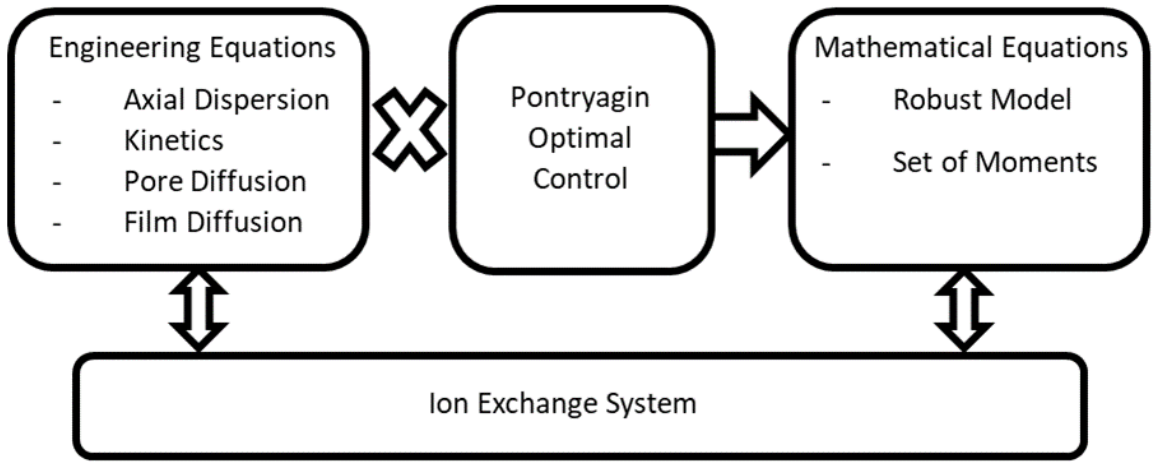

Figure 4. Adding optimal control on the ion exchange system via the use of moments rather than the various engineering equations.

Pontryagin's maximum principle employs the Hamiltonian by imposing the optimality conditions for the objective function (Artstein 2011). Optimal control is applied using a manageable system parameter such as the flow rate. With the absence of system uncertainties, derivatization of the Hamiltonian function from equation (12) with the flow rate (Q) can be represented with its two coordinates, the state variables, and the adjoint variables, using the following chain rule expression (Kirti M. Yenkie and Diwekar 2012):

$$\frac{dH}{d\chi} = \frac{dH}{dQ} = \sum_{i=1}^{n}\frac{dH}{dy_i}\frac{dy_i}{dQ} + \sum_{i=1}^{n}\frac{dH}{dz_i}\frac{dz_i}{dQ} = -\sum_{i=1}^{n}\frac{dz_i}{dt}\theta_i + \sum_{i=1}^{n}\frac{dy_i}{dt}\phi_i \qquad (25)$$

where $\phi_i$ & $\theta_i$ are new variables used to represent $dz_i/dQ$ & $dy_i/dQ$ respectively (Benavides and Diwekar 2012), and $dH/dy_i$ & $dH/dz_i$ are defined as $-dz_i/dt$ & $dy_i/dt$ respectively (Chiang 1992) as mentioned previously. Optimality is achieved when equation (25) is minimized confirming that the flow rate adjustment achieved the desired optimized objective. Therefore, an optimal flow rate is computed for every time increment ($\Delta t$) using equation (25).

While the state variables, in equations (13) to (16), are integrated forward from the initial time to the final time with the initial values being $y_i(0) \in [y_1(0)\ 0\ 0\ 0]$, the adjoint variables, in equations

(18) to (21), are integrated backward from the final time to the initial time. To find the terminal values of the adjoint variables, the derivative of the objective function in equation (24) is calculated and equated to zero (confirming that the maximum value was achieved):

$$\frac{d(J_{max})}{dt} = -C_0 Q_{final}\left(y_2(Q_{final}) - \sqrt{y_3(Q_{final})}\right)\frac{dy_1}{dt} + C_0 Q_{final}\left(1 - y_1(Q_{final})\right)\frac{dy_2}{dt} -$$

$$0.5\frac{C_0 Q_{final}}{\sqrt{y_3}}\frac{dy_3}{dt} = a\frac{dy_1}{dt} + b\frac{dy_2}{dt} + c\frac{dy_3}{dt} = 0 \qquad (26)$$

Therefore, the final values for the adjoint variable are extracted from equation (26) as $z_i$ ($t_{final}$) ∈ [$a$ $b$ $c$ 0]. Like the adjoint variables and the state variables, $\phi_i$ (representing $dz_i/dQ$) is integrated backward, from $\phi_i(t_{final})$ ∈ [0 0 0 0], and $\theta_i$ (representing $dy_i/dQ$) is integrated forward, from $\theta_i(0)$ ∈ [0 0 0 0].

If either the state variables, the adjoint variables, or the new variables, $\phi_i$ and $\theta_i$, are difficult to integrate, the Euler method of integration is used for each time step as seen below.

$$y_i(t+1) = y_i(t) + \Delta t \frac{dy_i}{dt} \quad \text{and} \quad \theta_i(t+1) = \theta_i(t) + \Delta t \frac{d\theta_i}{dt} \qquad (27)$$

$$z_i(t-1) = z_i(t) - \Delta t \frac{dz_i}{dt} \quad \text{and} \quad \phi_i(t-1) = \phi_i(t) - \Delta t \frac{d\phi_i}{dt} \qquad (28)$$

where $\Delta t$ – time increment used in the integration of the Euler method ($\Delta t < \frac{t_{final}}{500}$ for small steps). From equations (27) and (28), the pathways for both the state variables, $y_i$, and the adjoint variables, $z_i$, are calculated from initial values to final values and from final values to initial values respectively. The flow rate will be changed until changes in the Hamiltonian with the flow rate, from equation (25), are annulled. Therefore, a restricted pathway for an optimum flow rate is created to achieve such an objective.

**2.6. Capturing the uncertainties with the Ito process:** Up to this point, we approached our process as an absolute measurement with high accuracy. We did not consider any process fluctuations or uncertainties that exist in real-life situations. Fluctuations can come from many sources, such as concentration and other system measurements. Such uncertainties can drastically influence the ion exchange process duration and can even lead to unexpected leakage of chromates into the effluent stream. These fluctuations can make the optimal control strategy more challenging to implement. Therefore, it is essential to capture the uncertainties of the system conditions via the use of a stochastic process model. The approach to implementing the stochastic control strategy requires three separate steps.

i) The first step consists of evaluating the effect of parameter fluctuations on the process effluent. These fluctuations can be simulated via a probabilistic distribution with mean size and distribution.

ii) The second step applies a stochastic process model that best represents the results obtained from the first step. The Ito process model equation of Brownian motion with drift (Diwekar 2008) is used to display the effluent concentration fluctuations as a result of the process variations.

iii) The third step modifies the deterministic control previously shown in sections 2.3 to 2.5 by using additional parameters contained in the Ito stochastic process model.

Therefore, as a first step, we consider the variation of the inlet concentrations of $\Delta C_0 = \pm 10\%$ and the fluctuations of the model parameters of $\Delta q_m$ and $\Delta K_T = \pm 30\%$ on the ion exchange process (Ghanem, Jerpoth, and Yenkie 2022). The variations of all four moments will be plotted to mark the maximum and minimum of at least 100 variations. The average of the fluctuations and the results of the optimal deterministic plots are also added. The goal is to use a stochastic equation that will cover the deterministic plot and the average of the fluctuations as presented in the results section.

Stochastic processes are equations with parameters that can change randomly with time (Kao 2019). One of the most applied stochastic processes is the random walk process, also called the Wiener process (Szabados 2010), where one begins at a known initial value and takes equiprobable steps in any direction. Therefore, the next step becomes the new initial step when further equiprobable steps are taken, a critical Markov property (Briskot et al. 2019) where all future steps depend on the earlier ones. Such a stochastic process results in a normally distributed probability of directional steps with variance and a mean. The Brownian motion process (applying probability to the Markov chain) encompasses such an approach (Brereton 2014).

One example within the Brownian motion method is the Ito process reflecting the Brownian motion with drift (Kirti M. Yenkie and Diwekar 2018) as defined below.

$$dx = \alpha \, dt + \kappa \, dz \tag{29}$$

where dx – shift in the function, $\alpha$ – drift parameter, dt – shift in the variable, $\kappa$ – deviation parameter, and dz – increment in the Wiener process (reflecting the step changes for variations or system fluctuations) as seen below.

$$dz = \varepsilon \, (\Delta t)^{1/2} \tag{30}$$

where $\varepsilon$ – random number with a normal distribution of zero as the mean and a standard deviation of 1, and $\Delta t$ – increment of time that is used when creating a stochastic approximation of the discrete values as seen in equation (30) (K.M. Yenkie, Diwekar, and Linninger 2016).

$$x_t = x_{t-1} + \alpha \, \Delta t + \kappa \, \varepsilon \, (\Delta t)^{1/2} \tag{31}$$

The Ito process has been used successfully in reflecting fluctuations in other modes of purification such as crystallization (Kirti M. Yenkie and Diwekar 2012) and distillation (Lei, Li, and Chen 2003) where the stochastic equations were able to capture the uncertainties via similar behavior as the average of the system fluctuations. For this paper, equation (30) can be rewritten to reflect the moments, the drift parameter $\alpha = F = dy/dt$, and the deviation parameter $\kappa = g = \sqrt{\frac{\text{var}(\Delta y)}{\Delta t}}$ (Benavides and Diwekar 2012). When applying the Ito process in equation (30) to the moments, the following equations represent the four moments as they capture the same uncertainties of the system deviations. While the uncertainties for each moment can be captured with a different random value for ε, we utilized a standard normal distribution of mean zero and a distribution of 1 for simplicity.

$$\mu_0(t) = \mu_0(t-1) + F_{\mu_0} \Delta t + g_{\mu_0} \varepsilon (\Delta t)^{1/2}$$
$$= \text{captures the fluctuation in the Thomas model} \qquad (32)$$

$$\mu_1(t) = \mu_1(t-1) + F_{\mu_1} \Delta t + g_{\mu_1} \varepsilon (\Delta t)^{1/2}$$
$$= \text{captures the fluctuation in the residence time of the peak} \qquad (33)$$

$$\mu_{2c}(t) = \mu_{2c}(t-1) + F_{\mu_{2c}} \Delta t + g_{\mu_{2c}} \varepsilon (\Delta t)^{1/2}$$
$$= \text{captures the fluctuation in the variance of the peak} \qquad (34)$$

$$\mu_{3c}(t) = \mu_{3c}(t-1) + F_{\mu_{3c}} \Delta t + g_{\mu_{3c}} \varepsilon (\Delta t)^{1/2}$$
$$= \text{captures the fluctuation in the skewness of the peak} \qquad (35)$$

If the variance is ignored (g = 0), these equations would represent a deterministic model. In case the Ito process of Brownian motion with drift does not cover the fluctuations well, a different stochastic process such as the Ito process with reverting mean (Shastri and Diwekar 2006), can be used in cases where fluctuations occur around a mean value, as shown in equation (35).

$$\mu_{3c}(t) - \mu_{3c}(t-1) = \eta \, \mu_{3c_{ave}} \Delta t - \eta \, \mu_{3c}(t-1) \Delta t + \sigma \varepsilon \sqrt{\Delta t}$$
$$= \eta \, (\mu_{3c_{ave}} - \mu_{3c}(t-1)) \, \Delta t + \sigma \varepsilon \sqrt{\Delta t} \qquad (36)$$

where $\mu_{3c_{ave}}$ – the average to which the third moment reverts and η – the speed of the reversion to the average value.

To simplify the representation of the reverting mean as a differential equation like the Brownian motion with drift, we performed the following calculations and assumptions to simplify our optimal control approach:

$$\mu_{3c_{ave}} - \mu_{3c}(t-1) = \frac{\mu_{3c_1} + \mu_{3c_2} + \mu_{3c_3} + \cdots + \mu_{3c_n}}{n} - \mu_{3c}(t-1) =$$

$$\frac{\left(\mu_{3c_1} - \mu_{3c}(t-1)\right) + \left(\mu_{3c_2} - \mu_{3c}(t-1)\right) + \left(\mu_{3c_3} - \mu_{3c}(t-1)\right) + \cdots + \left(\mu_{3c_n} - \mu_{3c}(t-1)\right)}{n} = \Delta\mu_{3c_{ave}} \quad (37)$$

where $\Delta\mu_{3c_{ave}}$ – the average of the difference between the reverting mean and a preceding value. Since η is defined as the speed of reversion to the mean, it can be approximated as the average speed to revert.

$$\eta = \frac{1}{\Delta t_{ave}} \quad (38)$$

Therefore, the following approximation can be taken when proceeding with the control strategy.

$$\eta \left(\mu_{3c_{ave}} - \mu_{3c}(t-1)\right) = \frac{\Delta\mu_{3c_{ave}}}{\Delta t_{ave}} \sim \frac{d\mu_{3c}}{dt} = F_{\mu_{3c}} \quad (39)$$

Equation (39) allows us to use equation (35) (Brownian motion with drift) when applying the stochastic process control in the next section.

## 2.7. Stochastic maximum principle

The objective of this application remains is to maximize the amount of chromate removed by the ion exchange system at an optimal time $t_m$ with minimal deviation σ via flow rate control. But we now must start capturing the uncertainties via the use of the stochastic equations introduced in section 2.6 while maximizing the objective function (Rodriguez-Gonzalez et al. 2019).

When working with stochastic processes, we apply Ito's lemma (Benavides and Diwekar 2013), the stochastic calculus counterpart of the chain rule to an ordinary differential equation. Therefore, Ito's lemma is like a Taylor series expansion of a derivative. As the objective function, J, is a function of the moments and time, dJ can be represented by the following expansion:

$$dJ = \frac{dJ}{dt}dt + \frac{dJ}{dx}dx + \frac{1}{2}\frac{d^2J}{d^2t}dt^2 + \frac{1}{2}\frac{d^2J}{d^2x}dx^2 + \frac{1}{6}\frac{d^3J}{d^3t}dt^3 + \frac{1}{6}\frac{d^3J}{d^3x}dx^3 + \cdots \quad (40)$$

where t – time, and x – moment.

When expanding equation (40) to the first order with time and second order with the moments, substituting equation (29) will lead to the following equation assuming $dtdz \sim 0$, $dt^2 \sim 0$, and $dz^2 \sim dt$.

$$dJ = \left(\frac{dJ}{dt} + \alpha\frac{dJ}{dx} + \frac{1}{2}\kappa^2\frac{d^2J}{d^2x}\right)dt + \kappa\frac{dJ}{dx}dz \qquad (41)$$

With the optimality condition shown by equation (26), equation (41) can be rewritten via Ito Lemma as follow:

$$0 = \frac{dJ}{dt} + \max_Q\left(\alpha\frac{dJ}{dx} + \frac{1}{2}\kappa^2\frac{d^2J}{d^2x}\right) \qquad (42)$$

Equation (42) can then be expanded into the various moments with $\alpha = F_{y_i}$ and $\kappa = g_{y_i}$. The newly expanded equation is seen below:

$$0 = \frac{dJ}{dt} + \max_Q\left(\sum_{i=1}^{4}\frac{dJ}{dy_i}F_{y_i} + \sum_{i=1}^{4}\frac{d^2J}{d^2y_i}\frac{g_{y_i}^2}{2}\right) \qquad (43)$$

Defining the adjoint variables $z_i = dJ/\partial d$ and $\omega_i = d^2J/d^2y_i$, the extended stochastic Hamiltonian function can be written as follow (Abbasi and Diwekar 2013):

$$H = \sum_{i=1}^{4}\left[z_iF_i + \omega_i\frac{g_{y_i}^2}{2}\right] \qquad (44)$$

where $g_{y_i} = \sqrt{\frac{var(\Delta y_i)}{\Delta t}} \qquad (45)$

The new adjoint variables are defined by their derivatives by the following equations (Gallego 2013) which are capable of capturing the uncertainties via the use of the Ito process of Brownian motion with drift:

$$\frac{dz_i}{dt} = -\sum_{k=1}^{4}\left[\frac{dF_k}{dy_i}z_k + \frac{1}{2}\frac{dg_k^2}{dy_i}\omega_k\right] \qquad (46)$$

$$\frac{d\omega_i}{dt} = -\sum_{k=1}^{4}\left[2\omega_k\frac{dF_k}{dy_i} + \frac{d^2F_k}{dy_i^2}z_k + \frac{1}{2}\frac{d^2g_k^2}{dy_i^2}\omega_k\right] \qquad (47)$$

The approach via stochastic optimization is similar to deterministic optimization except for the use of two adjoint variables, $z_i$ and $\omega_i$, to capture the absolute and the variance of the pathway towards the objective function.

**2.8. Running the optimization program**: Further details for the presented calculations in sections 2.1 to 2.7 can be found in the supplementary section S1. The pump flow rate is bounded between 0.42 liters/hr and 1.27 liters/hr. Considering that the inlet chromate concentration enters at 20 ppb, the flow rate, initially set at the slow flow of 0.42 liters/hr, is expected to change to optimize the removal of chromates.

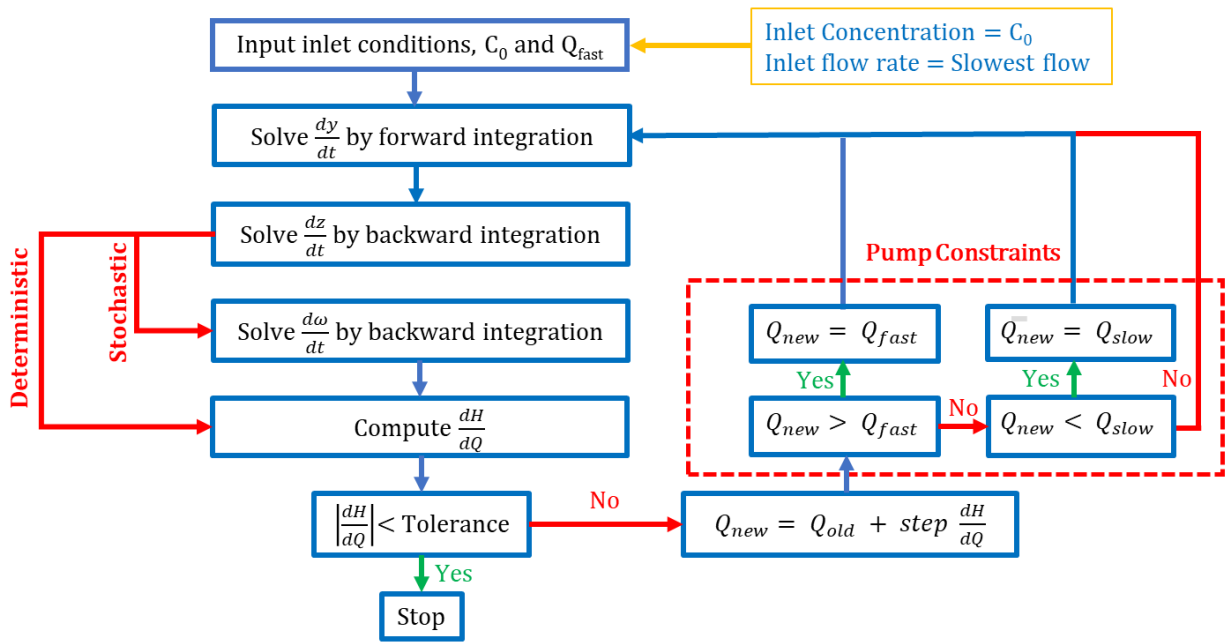

Figure 5. Flowchart to achieve deterministic and stochastic optimal control on adjusting the flow rate to maximize the performance of the chromate removal system.

Figure 5 summarizes the algorithm using Matlab as the optimization tool for the deterministic and stochastic approach. The Matlab optimization programs used can also be seen in the supplementary file, section S2. The difference between a stochastic approach and a deterministic approach is that the values of $\omega$ are ignored, corresponding to no uncertainties (stochastic step $dz = 0$). As seen in Figure 5, the decision to change the flow rate at each time step is done via a comparison of the Hamiltonian flow rate gradient with the tolerance (dH/dQ < Tolerance ~ $10^{-9}$). The decision is made to either stop the program as the Hamiltonian gradient met the tolerance level, or the program continues making changes to the flow rate and recalculating the state and adjoint variables keeping the flow rate bounded by its lower and upper value. It is important to note that the value of the step change should be small enough not to miss the optimality condition. The derivation of the dH/dQ value can be seen below for the stochastic approach that captures the state variable and the two adjoint variables (Kirti M. Yenkie and Diwekar 2012):

$$\frac{dH}{dQ} = \sum_{i=1}^{4} \frac{dH}{dy_i}\frac{dy_i}{dQ} + \sum_{i=1}^{4} \frac{dH}{dz_i}\frac{dz_i}{dQ} + \sum_{i=1}^{4} \frac{dH}{d\omega_i}\frac{d\omega_i}{dQ}$$

$$= \sum_{i=1}^{4} \frac{dH}{dy_i}\theta_i + \sum_{i=1}^{4} \frac{dH}{dz_i}\phi_i + \sum_{i=1}^{4} \frac{dH}{d\omega_i}\Omega_i \tag{48}$$

where $\theta_i$, $\phi_i$, and $\Omega_i$ are new functions to simplify the derivatization of $y_i$, $z_i$, $\omega_i$ in respect to the flow rate Q. Calculations of $\theta_i$, $\phi_i$, and $\Omega_i$ are made via their corresponding derivatives seen below:

$$\frac{d\theta_i}{dt} = \frac{d\left(\frac{dy_i}{dQ}\right)}{dt} = \frac{d\left(\frac{dy_i}{dt}\right)}{dQ} = \frac{dF_{y_i}}{dQ} \tag{49}$$

$$\frac{d\phi_i}{dt} = \frac{d\left(\frac{dz_i}{dQ}\right)}{dt} = \frac{d\left(\frac{dz_i}{dt}\right)}{dQ} \tag{50}$$

$$\frac{d\Omega_i}{dt} = \frac{d\left(\frac{d\omega_i}{dQ}\right)}{dt} = \frac{d\left(\frac{d\omega_i}{dt}\right)}{dQ} \tag{51}$$

As previously mentioned, the state variables, $y_i$ and $\theta_i$, are integrated forward from the initial conditions, [$y_1(0)$ 0 0 0] and [0 0 0 0] respectively. The adjoint variables, $z_i$ and $\omega_i$ as well as $\phi_i$ and $\Omega_i$, are integrated backward from the final conditions stated below based on the final optimum conditions of equation (26).

$$z(t_f) = [a, b, c, 0] =$$

$$\left[-\left(\mu_1(t_f) - \mu_{2c}(t_f)^{\frac{1}{2}}\right)C_0Q_f, \left(1 - \mu_0(t_f)\right)C_0Q_f, -\frac{(1-\mu_0(t_f))}{2\,\mu_{2c}(t_f)^{\frac{1}{2}}}C_0Q_f, 0\right] \quad (52)$$

$$\omega(t_f) = [0\ 0\ 0\ 0],\ \phi(t_f) = [0\ 0\ 0\ 0],\ \text{and}\ \Omega(t_f) = [0\ 0\ 0\ 0] \quad (53)$$

Similarly to the deterministic approach, the integrations were done using the Euler steps (Kirti M. Yenkie and Diwekar 2012) seen below:

$$y_{t+1} = y_t + \Delta t\, F_y \quad \text{and} \quad \theta_{t+1} = \theta_t + \Delta t\, \frac{d\theta}{dt} \quad (54)$$

$$z_{t-1} = z_t - \Delta t\, \frac{dz}{dt} \quad \text{and} \quad \phi_{t-1} = \phi_t - \Delta t\, \frac{d\phi}{dt} \quad (55)$$

$$\omega_{t-1} = \omega_t - \Delta t\, \frac{d\omega}{dt} \quad \text{and} \quad \Omega_{t-1} = \Omega_t - \Delta t\, \frac{d\Omega}{dt} \quad (56)$$

where $\Delta t$ – time step for the Euler method integration.

### 3. Results and Discussions:

The results section is split into three sections, 3.1 on the deterministic optimization results, 3.2 on capturing uncertainties, 3.3 on implementing the stochastic optimization, and 3.4 on the comparison between the two optimization approaches.

**3.1. Deterministic optimization results:** When we ran the solver program within Matlab (optimization program seen in supplementary section S2), we used more than 24000 iterations to achieve the optimal conditions. The results of such iterations are displayed in Figure 6. The iterations stopped once the changes in the Hamiltonian with the flow rate fell below the tolerance level (picked at $10^{-9}$). The program adjusted the flow rate periodically (every time step $\Delta t$) to maximize the objective function as previously described in Figure 5.

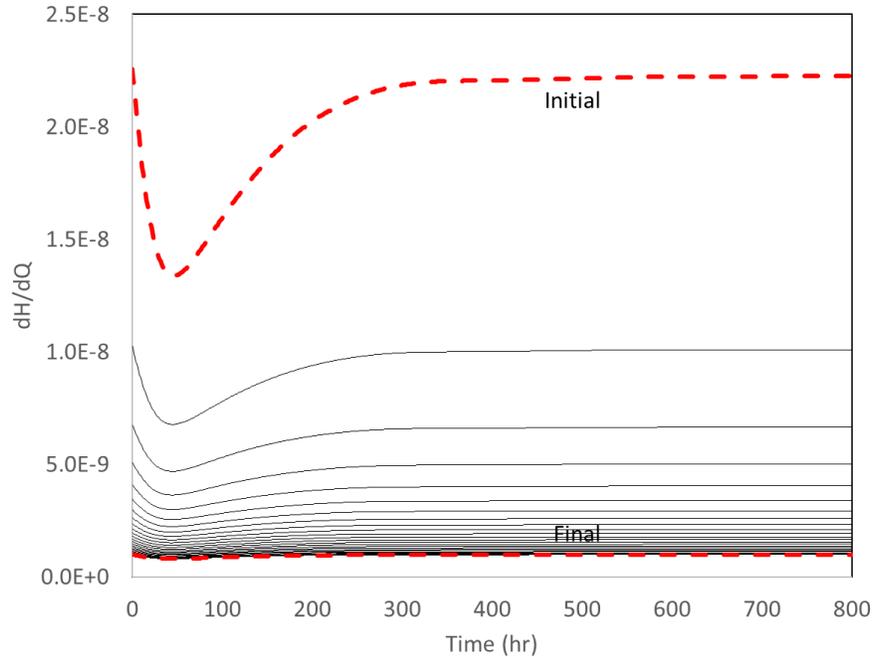

Figure 6. After 24040 iterations of the Hamiltonian differential with time, the optimal condition was achieved with an inlet chromate concentration of 20 ppb and an initial flow rate of 1.27 l/hr.

From Figure 6, the system equilibrated after 200 hours when the initial conditions consisted of an inlet chromate concentration of 20 ppb and a starting flow rate of 0.42 l/hr. That equilibration demonstrates that optimum conditions have been reached where flow rate changes are no longer necessary.

From Figure 7, both the Hamiltonian and the objective functions for the deterministic approach reached their optimum value around 200 hours. The Hamiltonian reached its constant minimum value, confirming that system optimality was reached, while the objective function peaked at its maximum value. That similarity demonstrates the effectiveness of the Hamiltonian to capture the optimality of the objective function.

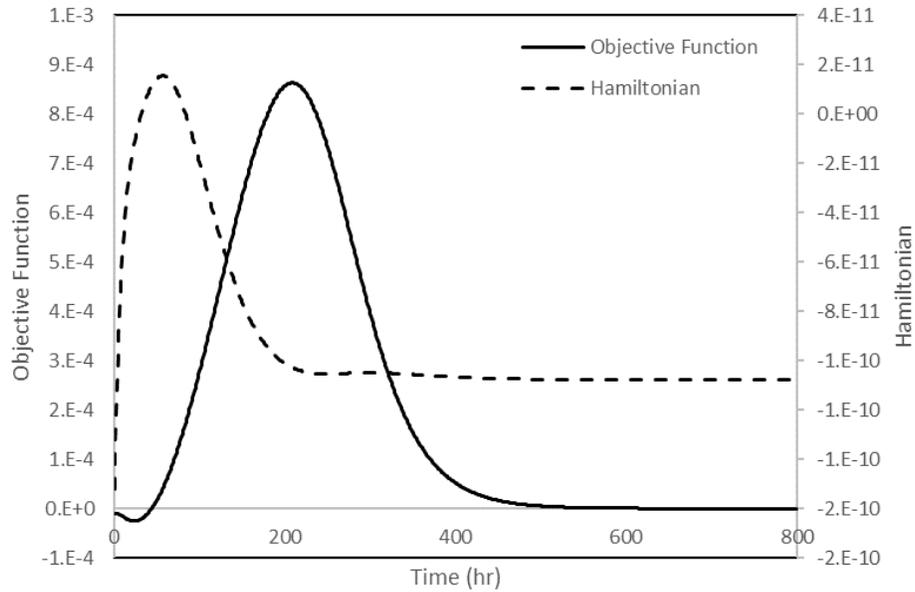

Figure 7. The system reached its maximum objective and minimum Hamiltonian after 200 hours of deterministic optimal control.

Around the 200 hours mark, there is no longer a need to change the flow rate since the final objective, $J_{max}$, was reached. The optimized dynamic flow rate path until the objective function is reached is seen in Figure 8.

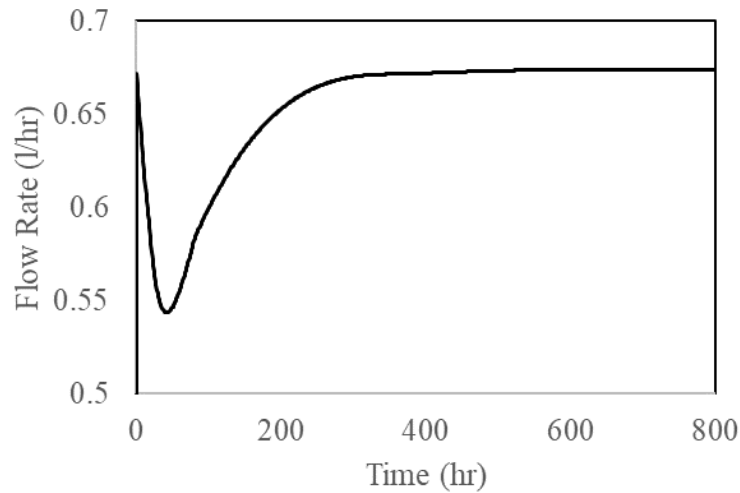

Figure 8. The optimal flow starts with an inlet chromate concentration of 20 ppb at the slowest flow rate of 0.42 l/hr.

From Figure 8, the final deterministic optimal flow rate of 0.67 l/hr was reached after 200 hours. The initial optimal flow rate increased to about 0.67 l/hr and then started slowing down until reaching a flow rate of 0.54 l/hr after about 50 hours before shifting upwards to maximize the ion exchange resin performance. The effect of slowing down the flow tends to increase the $t_m$ or first moment, while the effect of speeding up the flow tends to decrease the variance or the second moment, an essential characteristic of the objective function $(t_m - \sigma)$. The initial increase in the flow rate is because our objective function allows an acceptable amount of chromate to leak out before flow reduction becomes necessary. Our objective function also seeks the maximum use of the ion exchange resin via processing more chromate feed by keeping a high flow rate for as long as possible. But the flow is increased at the end as the system is seeking a quick failure while maintaining a low chromate leakage (the system is getting as much chromate extracted) due to the negative sign of the standard deviation in the objective function. The effect of the optimal flow rate is seen in Figure 9, demonstrating that optimality is achieved between the two flow rate boundaries of 0.42 l/hr and 1.27 l/hr.

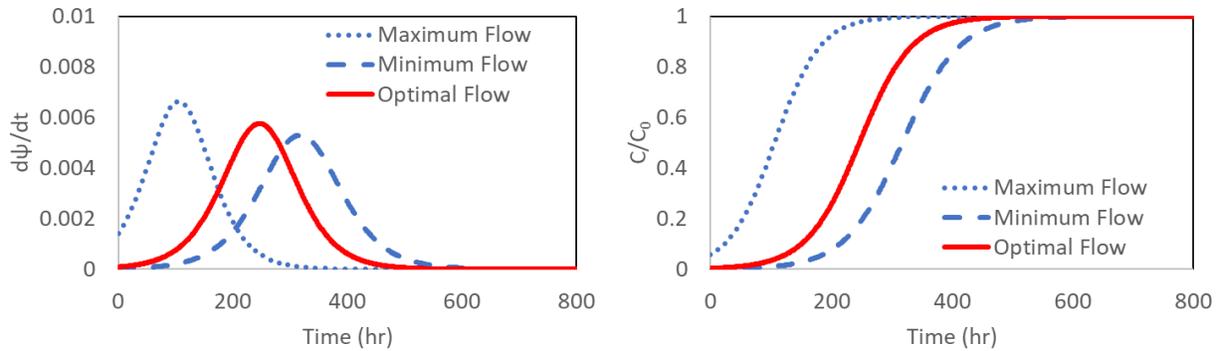

Figure 9. Effect of optimizing flow rates to maximize chromate removal between 2 flow rate boundaries.

To understand the sensitivity of the purification optimal control to the initial conditions, we varied the inlet chromate conditions and the initial flow rates by 10%, as seen in Figures 10 and 11.

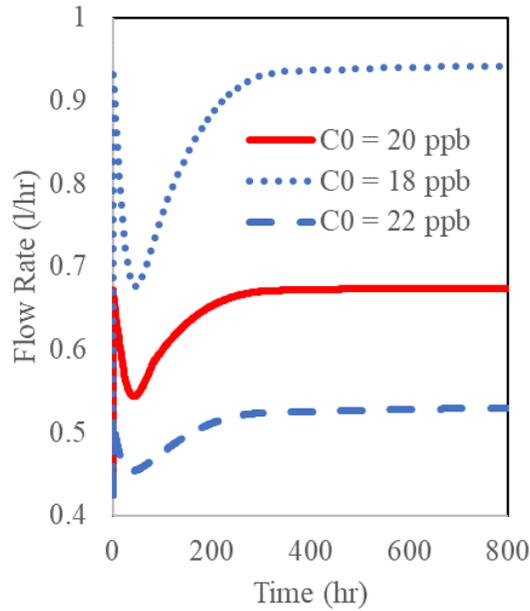

Figure 10. Effect of deviating the inlet chromate concentration by ±10% on the optimal conditions of the system.

Figure 10 demonstrates that increases in the inlet chromate concentration require the flow rate to run at lower flows to avoid early breakthroughs. The reverse is also expected when decreasing the inlet chromate concentration allowing the system to run at a higher flow rate since the risk for early breakthrough diminishes.

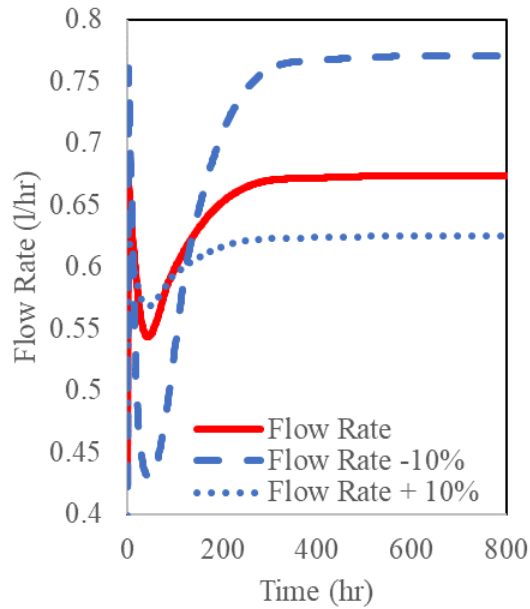

Figure 11. Effect of deviating the initial flow rate by ±10% on the optimal conditions of the system.

Figure 10 also demonstrates that the sensitivity of the system decreases as the inlet chromate concentration increases, forcing the optimal flow rate to stay close to its lowest limit.

Figure 11 confirms that higher initial flow rates will experience less significant flow rate modifications than lower initial flow rates. The reason for the larger flow rate modification with a lower initial flow rate is due to our objective function dependence on the first moment and second moment. At the specific time of $t_m - \sigma = 211$ hours, $C/C_0 = 14\%$. Therefore, when the concentration ratio is below 14%, the flow rate can run faster (but slower when $C/C_0$ is expected to be above 14%).

Table 2 summarizes how predictive optimal control can maximize the volume of processed water and chromate removal. By applying optimal control at a 14% concentration ratio (corresponding to the optimal time of $t_m - \sigma$), larger processed volumes and a higher amount of extracted chromate, compared to the slow or fast flow rate, was achieved. If the process was run to its optimal time of 211 hours, varying the flow rate would achieve even a greater difference than running at the slower flow rate. Therefore, optimizing the flow rate at each time step helped to maximize the removal of such a hazardous chemical.

Table 2. Amount of feed processed at various breakthrough points using fast flow, slow flow, and optimal flow

| Process Conditions | Fast Flow 1.27 l/hr | Slow Flow 0.42 l/hr | Optimal Flow varying with time |
|---|---|---|---|
| Time to reach $C/C_0 = 14\%$ (in hr) | 38 | 233 | 211 |
| Volume processed to reach $C/C_0 = 14\%$ (in liters) | 48 | 98 | 128 |
| Chromate removed by Purolite A600 to reach $C/C_0 = 14\%$ (in g per liter of resin) | 0.084 | 0.180 | 0.234 |
| Volume processed to optimal time of 211 hours (in liters) | 48* | 89 | 128 |
| Chromate removed by Purolite A600 to the optimal time of 211 hours (in g per liter of resin) | 0.081* | 0.164 | 0.234 |

* Value was taken at 38 hours since the system achieved the maximum 14% $C/C_0$ allowed

**3.2. Results of capturing uncertainties:** One hundred random variations were plotted for all four moments, as seen in Figure 12. We considered variations of the inlet concentrations of $\pm$ 10% and fluctuations of the model parameters of $\pm$ 30% on the ion exchange process (Ghanem, Jerpoth, and Yenkie 2022). The average of all 4 moments and the optimal results for the deterministic approach seem to follow a similar behavior.

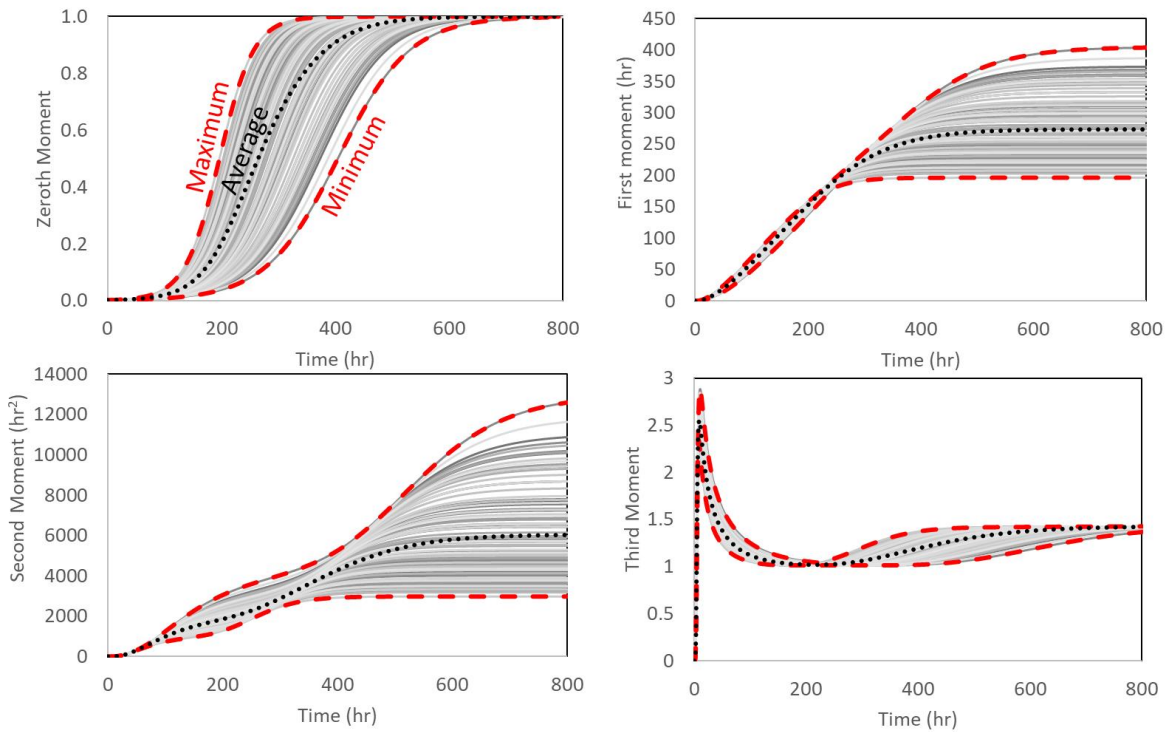

Figure 12. One hundred random variations of the inlet concentrations: $\Delta C_0 = \pm$ 10% and the Thomas model parameters: $\Delta q_m$ and $\Delta K_T = \pm$ 30% for all 4 moments.

Figure 13 demonstrates the four moments capturing time-dependent uncertainties with normalized distribution with zero mean and a standard deviation of one. The comparison of Figures 12 and 13 demonstrates that the Ito process can capture the uncertainties of the data using a randomness of zero mean and a distribution of 2. However, for the third moment, the use of the Ito process of Brownian motion with drift displayed more turbulence than expected from the ion exchange system since the sigmoidal function does not allow for too much skewness to occur (meaning that that it variates slightly around a mean). Since the Brownian motion with drift did not represent the third-moment fluctuations, the Ito process with reverting mean was plotted using an η =1.9 from equation (36) as seen in Figure 13. The reverting mean is a better stochastic representation of the third moment. Now that we captured the uncertainties using the Ito stochastic process, we can proceed to implement the stochastic optimization work.

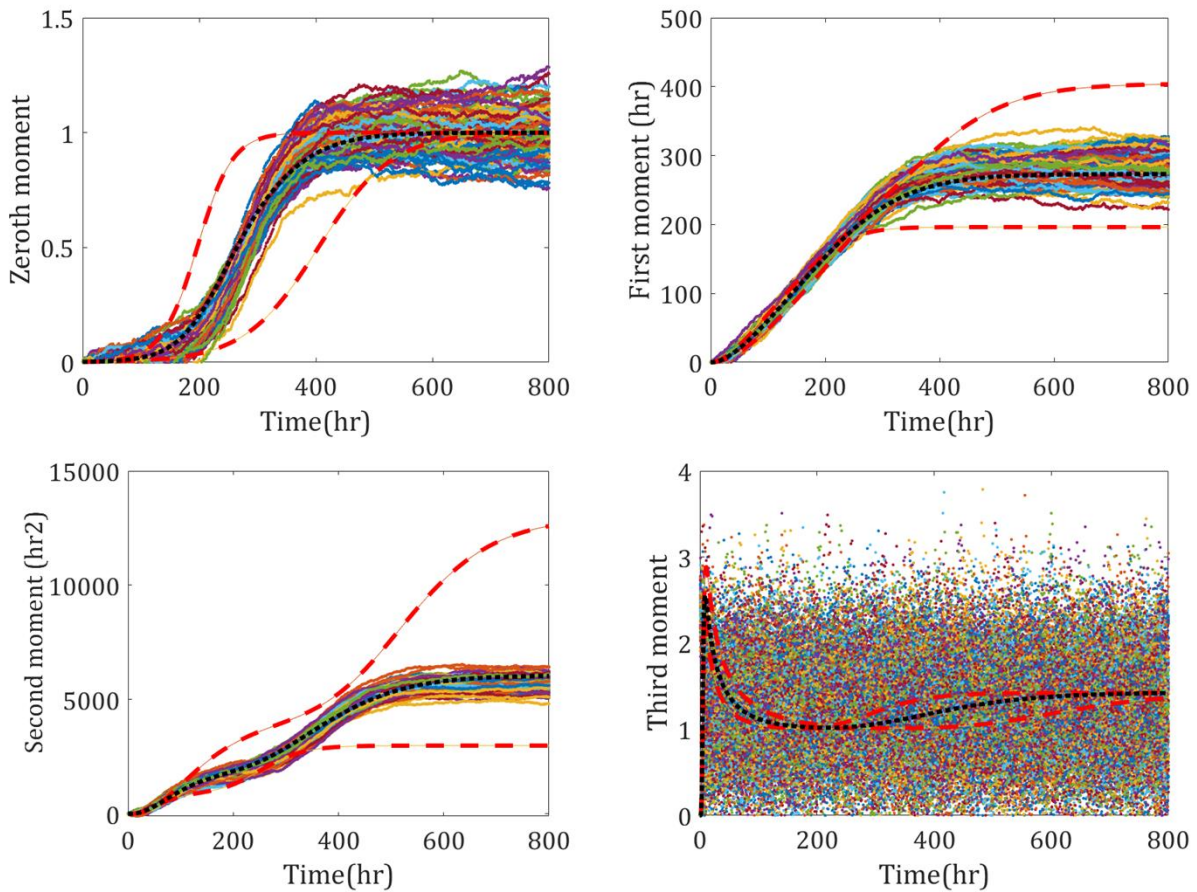

Figure 13. Ito process capturing the uncertainties with the four moments. The maximum, minimum, and average of the 100 variations of Figure 12 are included by comparison. Zeroth,

first, and second moment are used with the Brownian motion with drift while the third moment is used with the reverting mean.

**3.3. Stochastic optimization results:** When the Matlab program, summarized by the algorithm in Figure 5, was executed, the flow rate was adjusted constantly for thousands of iterations, minimizing the change of the Hamiltonian function with the flow rate at each time step (similar to what was done in the deterministic section). These results are seen in Figure 14.

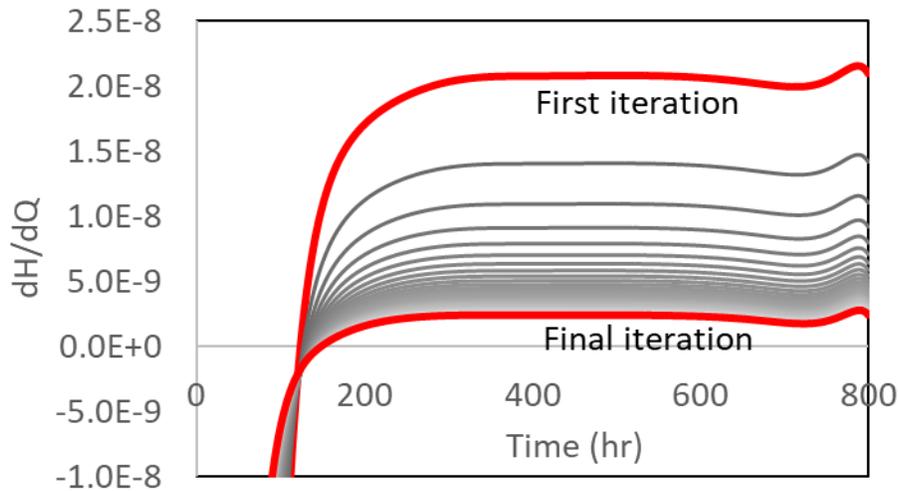

Figure 14. Optimization of the Hamiltonian changes with flow rates to achieve the objective function. Each line in the figure displays the results after 1000 iterations with the final line representing the results after 24000 iterations.

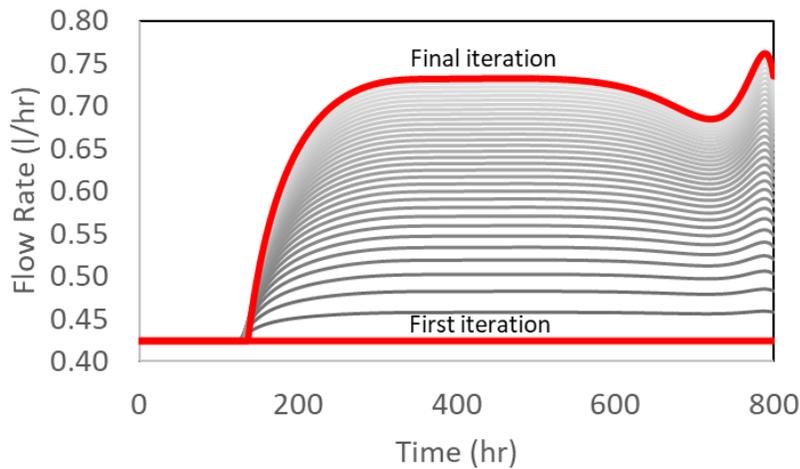

Figure 15. Changes in flow rates with 1000s of iterations with each line representing the flow rate every 1000 iterations.

From Figure 15, the final optimal flow rate of 0.72 l/hr was reached after 200 hours. The initial optimal flow rate remained at the initial 0.42 l/hr for about 150 hours before quickly increasing to the final flow of 0.72 l/hr, a considerable difference from what was previously achieved with the deterministic optimization. The stochastic approach considered all fluctuations including the impactful initial chromate concentration to control the flow rate. For example, if the initial chromate concentration was 10% higher, the flow rate would have decreased significantly, to its lower boundary, as demonstrated in Figure 15. Adding the effect of the 30% variation of the Thomas model parameters, the control system needed to run at a minimal flow rate to compensate for such fluctuations. The flow rate increases after 150 hours of the slow flow rate hold to enable as much chromate to be processed before unacceptable leakage occurs. It achieved its concentration ratio of 14% after 225 hours (as opposed to 211 hours for the deterministic case) since it applied the slow flow rate for an extended time as opposed to the deterministic approach.

**3.4. Comparison of the stochastic optimization with the deterministic optimization:** Figure 16 compares the objective function results for both the deterministic and the stochastic approaches. While the optimized flow rates were drastically different, the objective function for both approaches were similar in behavior, confirming that the paths can be different to achieve the same objective. The objective for the stochastic process was reached at 225 hours while the objective for the deterministic process was reached at 211 hours as noted in the previous section.

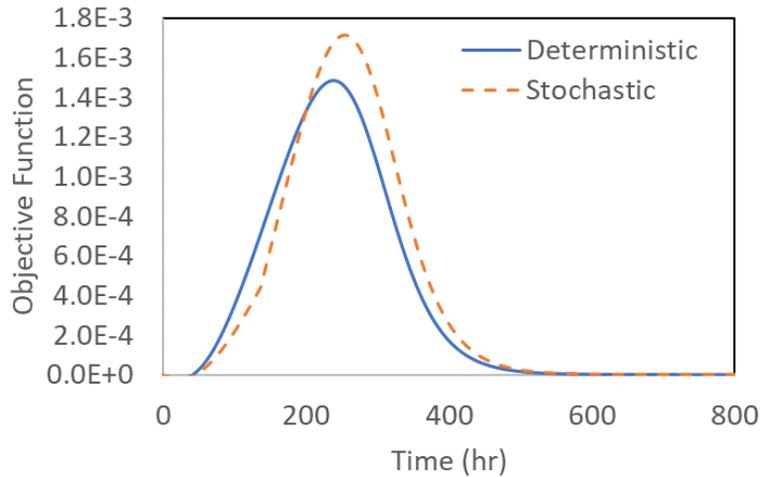

Figure 16. Comparison of the objective function for the deterministic and the stochastic approaches showing similarities.

While the objective functions are similar and the flow rates are different for both approaches, the effect on the concentration ratio can be seen in Figure 17, demonstrating a slightly longer run for the stochastic process to achieve the 14% concentration ratio.

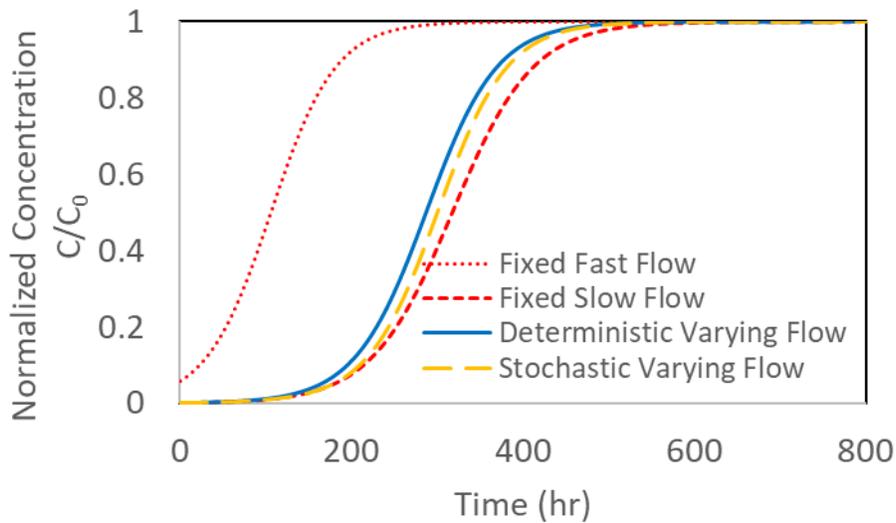

Figure 17. Impact of using a stochastic optimization over a deterministic optimization on the ion exchange system.

Since the results of the flow rate optimization were different between the stochastic and the deterministic approach, as summarized in Figure 18, we ran one hundred separate deterministic optimization programs to try capturing the uncertainty of the inlet concentrations ($\Delta C_0 = \pm\ 10\%$) and the model parameters ($\Delta q_m$ and $\Delta K_T = \pm\ 30\%$) and displayed them in Figure 20.

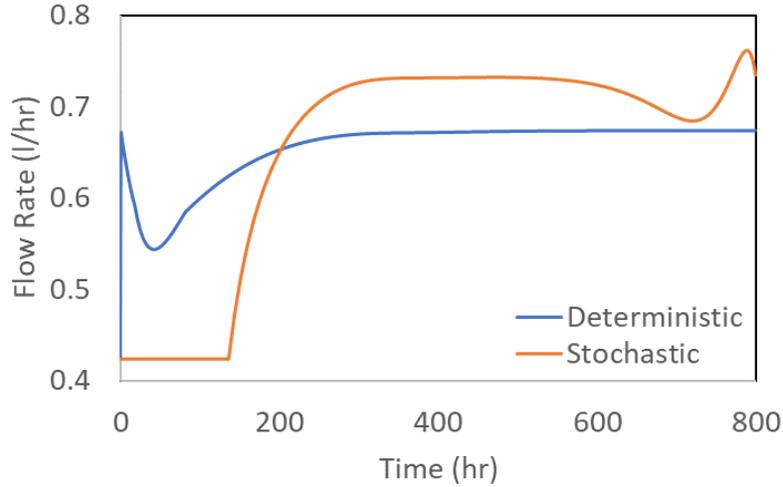

Figure 18. Comparison of the optimal flow rate between an absolute deterministic approach and a stochastic approach with a starting flow rate of 0.42 l/hr and an initial chromate concentration of 20 ppb.

Figure 19 demonstrates the sensitivity of the system to capture each uncertainty. The average of 100 deterministic optimal flow rates was getting closer to the stochastic optimized flow rate which seems to behave as the optimized version of thousands of deterministic optimizations. It took hours to run these 100 deterministic optimization programs versus a few seconds to run the stochastic optimization program. The stochastic approach displays the average of thousands of deterministic optimization pathways that capture each variation of the system, demonstrating the power of such an approach.

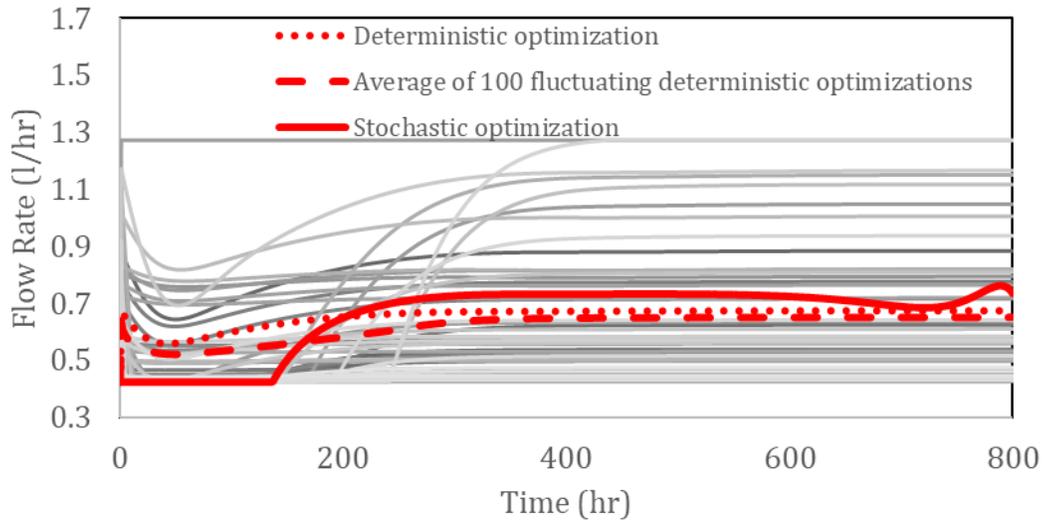

Figure 19. System flow rate behavior with each uncertainty comparing 100 different deterministic optimizations with the stochastic optimization approach.

With differences in flow rate control, a check on the value of the objective function was needed to confirm its effect on the system.

The results of both approaches are compared in Table 3 below. At the maximum objective function, the concentration ratio is calculated to be about 14% for both the deterministic and stochastic approaches. That 14% value corresponded to the dynamic binding capacity achieved at the optimal end time for both approaches.

Table 3. Comparison of the deterministic and stochastic optimal control on processing chromate solutions of 20 ppb Chromate and initial flow rate of 0.42 l/hr

| Flow rate to achieve $C/C_0 = 14\%$ | Volumes of chromate feed processed (liters) | Mass of chromates removed (grams of Cr per liter of resin) |
|---|---|---|
| Fixed Fast flow | 48 | 0.081 |
| Slow Flow Rate | 98 | 0.180 |
| Deterministic optimal flow | 128 | 0.234 |
| Stochastic optimal flow | 109 | 0.199 |

Table 3 confirms that optimal control can process more chromate solutions and remove more chromates than a fixed flow rate. This result demonstrates the necessity to vary the flow rate as a critical function to maximize the use of the ion exchange resin. But when comparing the results for optimal flow, the stochastic approach was more conservative as it captures the uncertainties of the system especially, when the inlet chromate concentration is higher or when the Thomas model parameters are more sensitive to fluctuations.

## 4. Conclusions:

Better control strategies are achieved when we measure inlet concentrations rather than outlet concentrations for hazardous compounds like chromates. For a successful implementation, we incorporate optimal control using a more accurate predictive model, such as the Thomas model, and Pontryagin's maximum principle via the use of moments. Our paper demonstrates that optimal chromate removal is achieved by running the system at varying flow rates rather than fixed flow to achieve more efficient use of our ion exchange resin. Optimal flow rate control resulted in processing more contaminated water and removing more chromate than what can be achieved with fixed flow rate settings.

Applying the stochastic approach demonstrated that the paths to optimal control can be different from the deterministic approach but still lead to higher chromate capture than the fixed flow rate. Furthermore, the amount of chromate removed via stochastic control was more conservative than the deterministic control since stochastic optimization captures the various uncertainties within the system.

With the successful implementation of the Hamiltonian via the method of moments, we demonstrate that smaller systems can also be used efficiently based on the inlet chromate concentrations. Since the ion exchange resin is rarely reused after concentrating the toxic chromates on its matrix, it is important to process more contaminated water optimally, resulting in less disposal of bulk ion exchange resins. We expect similar success with removing other toxic compounds such as arsenates, lead, and PFAS (perfluorooctanoic acid substances).

**Funding:**

This research did not receive any specific grant from funding agencies in the public, commercial, or not-for-profit sectors.

**Acknowledgments:**

The authors thank the Department of Chemical Engineering at Rowan University for their assistance in acquiring computational licenses and resources used in this study. Special thanks go to Swapana Jerpoth and Emmanuel Aboagye for their paper review and support as well as to Purolite for supporting my academic career.

**Supplementary Materials:**

A supplementary section was included with this paper to show the step-by-step calculations of the moments and Hamiltonian (S1) and the Matlab programs used for the optimization work (S2).

**References:**

Abbasi, Sheraz, and Urmila M Diwekar. 2013. "Characterization and Stochastic Modeling of Uncertainties in the Biodiesel Production." *Clean Technologies and Environmental Policy* 16: 79–94.

Artstein, Zvi. 2011. "Pontryagin Maximum Principle Revisited with Feedbacks." *European Journal of Control* 1: 46–54.

Balan, Catalin, Irina Volf, and Doina Bilba. 2013. "Chromium (VI) Removal from Aqueous Solutions by Purolite Base Anion-Exchange Resins with Gel Structure." *Chemical Industry and Chemical Engineering Quarterly* 19(4): 615–28.

Benavides, Pahola T., and Urmila Diwekar. 2012. "Optimal Control of Biodiesel Production in a Batch Reactor." *Fuel* 94: 218–26.

Benavides, Pahola T, and Urmila Diwekar. 2013. "Studying Various Optimal Control Problems in Biodiesel Production in a Batch Reactor under Uncertainty." *Fuel* 103: 585–92.

Biswas, Swarup, and Umesh Mishra. 2015. "Continuous Fixed-Bed Column Study and Adsorption Modeling: Removal of Lead Ion from Aqueous Solution by Charcoal Originated from Chemical Carbonization of Rubber Wood Sawdust." *Journal of Chemistry* 2015: 1–9.

Boscain, Ugo, and Benedetto Piccoli. 2005. "A Short Introduction to Optimal Control." In *Contrôle Non Linéaire et Applications: Cours Donnés à l'école d'été Du Cimpa de l'Université de Tlemcen / Sari Tewfit.*, , 19–66.

Brereton, Dr Tim. 2014. In *Stochastic Simulation of Processes, Fields and Structures*, Ulm University: Institute of Stochastics, 108–21.

Briskot, Till et al. 2019. "Prediction Uncertainty Assessment of Chromatography Models Using Bayesian Inference." *Journal of Chromatography A* 1587: 101–10.

Brito F et al. 1997. "Equilibria of Chromate(VI) Species in Acid Medium and Ab Initio Studies of These Species,." *Polyhedron* 16(21): 3835–46.

Carta, G., and A. Jungbauer. 2010. "Effects of Dispersion and Adsorption Kinetics on Column Performance." In *Protein Chromatography*, John Wiley & Sons, Ltd, 237–76.

Charola, Samir, Rahul Yadav, Prasanta Das, and Subarna Maiti. 2018. "Fixed-Bed Adsorption of Reactive Orange 84 Dye onto Activated Carbon Prepared from Empty Cotton Flower Agro-Waste." *Sustainable Environment Research* 28(6): 298–308.

Chiang, Alpha C. 1992. "Optimum Control: The Maximum Principle." In *Elements of Dynamic Optimization*, McGraw-Hill, 167–77.

Corder, GD, and PL Lee. 1986. "Feedforward Control of a Wastewater Plant." *Water Research.* 20(3): 301–9.

Costa, Max, and Catherine B. Klein. 2006. "Toxicity and Carcinogenicity of Chromium Compounds in Humans." *Critical Reviews in Toxicology* 36(2): 155–63.

de Dardel, François, and Thomas V. Arden. 2008. "Ion Exchangers." In *Ullmann's Encyclopedia of Industrial Chemistry*, ed. Wiley-VCH Verlag GmbH & Co. KGaA. Weinheim, Germany: Wiley-VCH Verlag GmbH & Co. KGaA.

Diwekar, Urmila. 2008. "Optimal Control and Dynamic Optimization." In *Introduction to Applied Optimization*, Springer Optimization and Its Applications, Boston, MA: Springer US, 215–77.

Ferreira, M.G.S., M.L. Zheludkevich, and J. Tedim. 2011. "9 - Advanced Protective Coatings for Aeronautical Applications." In *Nanocoatings and Ultra-Thin Films*, eds. Abdel Salam Hamdy Makhlouf and Ion Tiginyanu. Woodhead Publishing, 235–79.

Fogler, H. Scott. 2016. "Residence Time Distributions of Chemical Reactors." In *Elements of Chemical Reaction Engineering*, Boston: Prentice Hall, 767–806.

Gallego, Pahola Thathiana Benavides. 2013. "Optimal Control of Batch Production of Biodiesel Fuel under Uncertainty." University of Illinois.

Ghanem, Fred, Swapana S. Jerpoth, and Kirti M. Yenkie. 2022. "Improved Models for Chromate Removal Using Ion Exchangers in Drinking Water Applications." *Journal of Environmental Engineering* 148(5). https://doi.org/10.1061/(ASCE)EE.1943-7870.0001997 (March 20, 2022).

Goltz, Mark N., and Paul V. Roberts. 1987. "Using the Method of Moments to Analyze Three-Dimensional Diffusion-Limited Solute Transport from Temporal and Spatial Perspectives." *Water Resources Research* 23(8): 1575–85.

Hamdaoui, Oualid. 2009. "Removal of Copper(II) from Aqueous Phase by Purolite C100-MB Cation Exchange Resin in Fixed Bed Columns: Modeling." *Journal of Hazardous Materials* 161(2–3): 737–46.

Hamill, Patrick. 2018. *Lagrangians and Hamiltonians*. 4th printing. Cambridge University Press.

Han, Yongming et al. 2023. "Novel Long Short-Term Memory Neural Network Considering Virtual Data Generation for Production Prediction and Energy Structure Optimization of Ethylene Production Processes." *Chemical Engineering Science* 267. https://doi.org/10.1016/j.ces.2022.118372 (April 16, 2023).

Harmand, Jerome, Claude Lobry, Alain Rapaport, and Tewfik Sari. 2019. 3 *Optimal Control in Bioprocesses: Pontryagin's Maximum Principle in Practice*. First Edition. Wiley.

Hutcheson, John. 2006. "Ultrapure Water: Systems for Microelectronics." *Filtration & Separation* 43(5): 22–25.

Jawitz, James W. 2004. "Moments of Truncated Continuous Univariate Distributions." *Advances in Water Resources* 27(3): 269–81.

Kabir, G. 2008. "Removal of Chromate in Trace Concentration Using Ion Exchange From Tannery Wastewater." *Int. J. Environ. Res.* (2(4)): 377–84.

Kalaruban, Mahatheva et al. 2016. "Removing Nitrate from Water Using Iron-Modified Dowex 21K XLT Ion Exchange Resin: Batch and Fluidised-Bed Adsorption Studies." *Separation and Purification Technology* 158: 62–70.

Kao, P.C. Kao. 2019. "Brownian Motion and Other Diffusion Processes." In *An Introduction to Stochastic Processes*, Dover, 373–420.

Karakurt, Sevtap, Erol Pehlivan, and Serdar Karakurt. 2019. "Removal of Carcinogenic Arsenic from Drinking Water By the Application of Ion Exchange Resins." *Oncogen* 2(1). https://doi.org/10.35702/onc.10005 (August 18, 2019).

Lei, Zhigang, Chengyue Li, and Biaohua Chen. 2003. "Extractive Distillation: A Review." *Separation & Purification Reviews* 32(2): 121–213.

LeVan, M Douglas, Giorgio Carta, and Carmen M Yon. 1999. "Section 16: Adsorption and Ion Exchange." In *Perry's Chemical Engineers Handbook*, McGraw-Hill.

Li, Xue et al. 2016. "Chromium Removal From Strong Base Anion Exchange Waste Brines." *Journal - American Water Works Association* 108: 247–55.

———. 2016. "Meeting the New California MCL for Hexavalent Chromium with Strong Base Anion Exchange Resin." *American Water Works Association* 108(9): 474–81.

Lin, S.H., and C.D. Kiang. 2003. "Chromic Acid Recovery from Waste Acid Solution by an Ion Exchange Process: Equilibrium and Column Ion Exchange Modeling." *Chemical Engineering Journal* 92(1–3): 193–99.

Luo, Jian, Olaf A. Cirpka, and Peter K. Kitanidis. 2006. "Temporal-Moment Matching for Truncated Breakthrough Curves for Step or Step-Pulse Injection." *Advances in Water Resources* 29(9): 1306–13.

Millar, Graeme J., Sara J. Couperthwaite, Mitchell de Bruyn, and Chun Wing Leung. 2015. "Ion Exchange Treatment of Saline Solutions Using Lanxess S108H Strong Acid Cation Resin." *Chemical Engineering Journal* 280: 525–35.

Mustafa, Yasmen A, and Shahlaa E Ebrahim. 2010. "Utilization of Thomas Model to Predict the Breakthrough Curves for Adsorption and Ion Exchange." *Journal of Engineering* 16(4): 6206–23.

Nagaki, M, R D Hughes, J Y N Lau, and Roger Williams. 1991. "Removal of Endotoxin and Cytokines by Adsorbents and the Effect of Plasma Protein Binding." *Int J Artif Organs* 14(1): 43–50.

Nur, T. et al. 2015. "Nitrate Removal Using Purolite A520E Ion Exchange Resin: Batch and Fixed-Bed Column Adsorption Modelling." *International Journal of Environmental Science and Technology* 12(4): 1311–20.

Purolite Website, 1. 2021. "Purolite A600E." https://www.purolite.com/product/a600e. Accessed 2021.

Recepoğlu, Yaşar K. et al. 2018. "Packed Bed Column Dynamic Study for Boron Removal from Geothermal Brine by a Chelating Fiber and Breakthrough Curve Analysis by Using Mathematical Models." *Desalination* 437: 1–6.

Rodriguez-Gonzalez, Pablo T, Vicente Rico-Ramirez, Ramiro Rico-Martinez, and Urmila M Diwekar. 2019. "A New Approach to Solving Stochastic Optimal Control Problems." *Mathematics* 7(1207): 13.

Sharma, Deepak. 2022. "Simulated Moving Bed Technology: Overview and Use in Biorefineries." In *Biorefineries - Selected Processes*, IntechOpen, 1–21.

Shastri, Y, and U Diwekar. 2006. "An Optimal Control and Options Theory Approach to Forecasting and Managing Sustainable Systems." In *International Congress on Environmental Modeling and Software*, Burlington, VT.

Staby, Arne et al. 2007. "Comparison of Chromatographic Ion-Exchange Resins." *Journal of Chromatography A* 1164(1–2): 82–94.

Szabados, Tamas. 2010. "An Elementary Introduction to the Wiener Process and Stochastic Integrals." *arXiv:1008.1510 [math]*. http://arxiv.org/abs/1008.1510 (February 21, 2021).

Yenkie, K. M., and U.M. Diwekar. 2012. "Stochastic Optimal Control of Seeded Batch Crystallizer Applying the Ito Process." *Industrial & Engineering Chemistry Research*: 108–22.

———. 2018. "The 'No Sampling Parameter Estimation (NSPE)' Algorithm for Stochastic Differential Equations." *Chemical Engineering Research and Design* 129: 376–83.

Yenkie, K.M., U.M. Diwekar, and A.A. Linninger. 2016. "Simulation-Free Estimation of Reaction Propensities in Cellular Reactions and Gene Signaling Networks." *Computers & Chemical Engineering* 87: 154–63.

Yu, C., A. W. Warrick, and M. H. Conklin. 1999. "A Moment Method for Analyzing Breakthrough Curves of Step Inputs." *Water Resources Research* 35(11): 3567–72.